\documentclass[12pt,leqno]{article}

\usepackage{amsthm}

\usepackage{amsmath}
\usepackage{amscd}
\usepackage{amsmath}
\usepackage{amssymb}
\usepackage{latexsym}
\makeatletter

\@addtoreset{equation}{section}
\makeatother

\newtheorem{thm}{Theorem}
\newtheorem{pr}[thm]{Proposition}

\newtheorem{lm}[thm]{Lemma}

\newtheorem{cn}[thm]{Conjecture}

\newcommand{\sm}{\raisebox{2.33pt}{~\rule{6.4pt}{1.3pt}~}}

\input xy \xyoption{all} \CompileMatrices

\begin{document}

\title{On the proper push-forward
of the characteristic cycle\\
of a constructible sheaf}
\author{Takeshi Saito}

\maketitle

\begin{abstract}
We study the compatibility
with proper push-forward
of the characteristic cycles 
of a constructible complex
on a smooth variety over
a perfect field.
\end{abstract}

The characteristic cycle of
a constructible complex on
a smooth scheme over
a perfect field 
is defined as a cycle on the cotangent
bundle \cite{CC}
supported on the singular support
\cite{Be}.
It is characterized by
the Milnor formula \cite[(5.15)]{CC}
for the vanishing cycles
defined for morphisms to curves.

We study the compatibility
with proper push-forward.
First, we formulate Conjecture \ref{cnf*} on
the compatibility with proper direct image.
We prove it in some cases, for example,
morphisms from surfaces to
curves under a mild assumption
in Theorem \ref{thm1}.
We briefly sketch the idea of proofs,
which use the global index formula
(\ref{eqif})
computing the Euler-Poincar\'e characteristic.
For the compatibility Theorem \ref{thm1},
it amounts to prove a conductor
formula (\ref{eqaf})
at each point of the curve.
By choosing a point
and killing ramification
at the other points using Epp's
theorem \cite{Epp},
we deduce the conductor formula (\ref{eqaf})
from the index formula (\ref{eqif}).
We give a characterization
of characteristic cycle
in terms of functorialities
at the end of the article.

The research was partially supported
by JSPS Grants-in-Aid 
for Scientific Research
(A) 26247002.

\medskip



For the definitions and
basic properties 
of the singular support
of a constructible complex
on a smooth scheme over
a perfect field,
we refer to \cite{Be}
and \cite{CC}.

Let $k$ be a perfect field and 
let $\Lambda$ be a finite field
of characteristic $\ell$ invertible in $k$.
We say that a complex ${\cal F}$ of 
$\Lambda$-modules on the \'etale
site of $X$ is constructible
if the cohomology sheaf 
${\cal H}^q{\cal F}$ is constructible
for every $q$ and vanishes
except for finitely many $q$.
More generally, if $\Lambda$
is a finite local ring of
residue field $\Lambda_0$
of characteristic $\ell$ invertible in $k$,
the singular support and
the characteristic cycle
of a complex ${\cal F}$ of $\Lambda$-modules
of finite tor-dimension
equals to those of
${\cal F}\otimes_\Lambda^L\Lambda_0$.

Let ${\cal F}$ be a constructible complex on 
a smooth scheme $X$ over $k$.
The singular support
$SS{\cal F}$ is defined in \cite{Be}
as a closed conical subset of
the cotangent bundle $T^*X$
By \cite[Theorem 1.3 (ii)]{Be},
every irreducible component $C_a$
of the singular support
$SS{\cal F}=C=\bigcup_aC_a$ 
is of dimension $n$.
The characteristic cycle
$CC{\cal F}=\sum_am_aC_a$
is defined as a linear combination
with ${\mathbf Z}$-coefficients
in \cite[Definition 5.10]{CC}.
It is characterized by the
Milnor formula
\begin{equation}
-\dim {\rm tot}
\phi_u({\cal F},f)=
(CC{\cal F},df)_{T^*U,u}
\label{eqMil}
\end{equation}
for morphisms $f\colon U\to Y$
to smooth curves $Y$
defined on an \'etale neighborhood $U$
of an isolated characteristic point $u$.
For more detail on the notation, 
we refer to \cite[Section 5.2]{CC}.

We say that a constructible
complex ${\cal F}$ is {\em locally 
constant} if every cohomology sheaf
${\cal H}^q{\cal F}$ is locally constant.
In this case,
we have
\begin{equation}
CC{\cal F}=
(-1)^n{\rm rank}\ {\cal F}\cdot [T^*_XX]
\end{equation}
where $T^*_XX$
denotes the $0$-section
and $n=\dim X$
by \cite[Lemma 5.11.1]{CC}.
Assume $\dim X=1$
and 
let $U\subset X$ be a dense open subset
where ${\cal F}$ is
locally constant.
For a closed point $x\in X$,
the Artin conductor
$a_x{\cal F}$
is defined by
\begin{equation}
a_x{\cal F}
=
{\rm rank}\ {\cal F}|_U
-
{\rm rank}\ {\cal F}_x
+
{\rm Sw}_x{\cal F}
\label{eqaxF}
\end{equation}
where 
${\rm Sw}_x{\cal F}$
denotes the alternating sum of the
Swan conductor at $x$.
Then, by \cite[Lemma 5.11.3]{CC}, we have
\begin{equation}
CC{\cal F}
=
-\Bigl({\rm rank}\ {\cal F}\cdot [T^*_XX]
+\sum_{x\in X\!\! \sm \!\! U}a_x{\cal F}\cdot
[T^*_xX]\Bigr)
\end{equation}
where $T^*_xX$
denotes the fiber.

To state the compatibility with
push-forward, we fix some terminology
and notations.
We say that a morphism
$f\colon X\to Y$ of
noetherian scheme is proper (resp.\! finite)
on a closed subset $Z\subset X$
if its restriction $Z\to Y$
is proper (resp.\! finite) with respect to a
closed subscheme structure of $Z\subset X$.

Let $h\colon W\to X$ and
$f\colon W\to Y$ be 
morphisms of smooth schemes over $k$.
Let $C\subset X$ be a closed subset
such that
$f$ is proper on $h^{-1}(C)$
and let $C'=f(h^{-1}(C))\subset Y$ be
the image of $C$ by the algebraic
correspondence $X\gets W\to Y$.
If $\dim W=\dim X-c$,
the intersection theory
defines the pull-back and
push-forward morphisms
\begin{equation}
\begin{CD}
CH_\bullet (C)
@>{h^!}>>
CH_{\bullet -c}(h^{-1}(C))
@>{f_*}>>
CH_{\bullet -c}(C').
\end{CD}
\label{eqCH}
\end{equation}
We call the composition
the morphism defined by 
the algebraic
correspondence $X\gets W\to Y$.
If every irreducible component
of $C$ is of dimension $n$
and
if every irreducible component
of $C'$ is of dimension $m=n-c$,
the morphism (\ref{eqCH})
defines
a morphism
$Z_n(C)\to Z_m(C')$
of free abelian groups of cycles.

Let $f\colon X\to Y$
be a morphism of smooth schemes
over a field $k$.
Assume that every irreducible
component of $X$ is
of dimension $n$
and that every irreducible
component of $Y$ is
of dimension $m$.
Let $C\subset T^*X$ be
a closed conical subset.
The intersection
$B=C\cap T^*_XX$ with the $0$-section
regarded as a closed subset of $X$
is called the base of $C$.

Assume that $f\colon X\to Y$
is proper on the base $B$.
Then, the morphism
$X\times_YT^*Y\to T^*Y$
induced by $f$
is proper on the inverse image
$df^{-1}(C)$ by the canonical
morphism
$df\colon X\times_YT^*Y\to T^*X$.
Let $f_\circ C\subset T^*Y$
denote the image of
$df^{-1}(C)$ by $X\times_YT^*Y\to T^*Y$.
Then $f_\circ C\subset T^*Y$
is a closed conical subset.

By applying the construction of (\ref{eqCH})
to the algebraic correspondence
$T^*X\gets X\times_YT^*Y\to T^*Y$
and $C\subset T^*X$,
we obtain a morphism
\begin{equation}
f_*\colon 
CH_n(C)\to CH_m(f_\circ C).
\label{eqf*}
\end{equation}
Further if every irreducible
component of 
$C\subset T^*X$ is of dimension $n$
and if every irreducible
component of 
$f_\circ C\subset T^*Y$ is of dimension $m$,
we obtain a morphism
\begin{equation}
f_*\colon 
Z_n(C)\to Z_m(f_\circ C).
\label{eqf*Z}
\end{equation}

Now, 
let $f\colon X\to Y$ 
be a morphism 
of smooth schemes over $k$
and assume that every irreducible
component of $Y$ is
of dimension $m$.
The base $B\subset X$ of the singular support
$C=SS{\cal F}\subset T^*X$
equals the support of $X$
by \cite[Lemma 2.1(i)]{Be}.
Then, the direct image
\begin{equation}
f_*CC{\cal F}
\in 
CH_m(f_\circ C)
\label{eqfCC}
\end{equation}
of the characteristic cycle $CC{\cal F}$
is defined by the algebraic correspondence
$T^*Y\gets X\times_YT^*Y\to T^*X$.
Further if every irreducible
component of $f_\circ C\subset T^*Y$
is of dimension $m$,
the direct image
\begin{equation}
f_*CC{\cal F}
\in 
Z_m(f_\circ C)
\label{eqfCCZ}
\end{equation}
is defined as a linear combination of cycles.

\begin{cn}\label{cnf*}
Let $f\colon X\to Y$
be a morphism 
of smooth schemes 
over a perfect field $k$.
Assume that
every irreducible
component of $X$ is
of dimension $n$ and that
every irreducible
component of $Y$ is
of dimension $m$.
Let ${\cal F}$ be a constructible complex on $X$
and $C=SS{\cal F}$ be the singular support.
Assume that $f$ is proper on
the support of ${\cal F}$.

{\rm 1.}
We have
\begin{equation}
CCRf_*{\cal F}=
f_*CC{\cal F}
\label{eqcnf}
\end{equation}
in $CH_m(f_\circ C)$.

{\rm 2.}
In particular,
if every irreducible
component of $f_\circ C\subset T^*Y$ is of 
dimension $m$,
we have
an equality {\rm (\ref{eqcnf})}
of cycles.
\end{cn}

If $Y={\rm Spec}\ k$ is a point
and $X$ is proper over $k$,
the equality {\rm (\ref{eqcnf})}
is nothing but the index formula
\begin{equation}
\chi(X_{\bar k},{\cal F})
=
(CC{\cal F},T^*_XX)_{T^*X}.
\label{eqif}
\end{equation}
This is proved in \cite[Theorem 7.13]{CC}
under the assumption that
$X$ is projective.
For a closed immersion $i\colon X\to P$
of smooth schemes over $k$,
Conjecture \ref{cnf*}
holds \cite[Lemma 5.13.2]{CC}.
Hence,  for a proper morphism
$g\colon P\to Y$
of smooth schemes over $k$,
Conjecture \ref{cnf*}
for ${\cal F}$ 
and $f=g\circ i\colon X\to Y$
is equivalent to that
for $i_*{\cal F}$ and $g\colon P\to Y$.
For the singular support,
an inclusion
$SSRf_*{\cal F}
\subset f_\circ SS{\cal F}$ is proved in
\cite[Theorem 1.4 (ii)]{Be}.

\begin{lm}\label{lmfini}
Assume that $f$ is
finite on the support of ${\cal F}$.
Then Conjecture {\rm \ref{cnf*}.2}
holds.
\end{lm}

\proof{
We may assume that $k$ is algebraically closed.
Since the characteristic cycle
is characterized by the Milnor formula,
it suffices to show that
$f_*CC{\cal F}$ satisfies
the Milnor formula (\ref{eqMil}) for 
$Rf_*{\cal F}$.

Let $Z\subset X$ denote the support of ${\cal F}$.
Let $V\to Y$ be an \'etale morphism
and $g\colon V\to T$ be a morphism
to a smooth curve $T$ with
isolated characteristic point $v\in V$
with respect to $f_\circ C$.
By replacing $Y$ by $V$,
we may assume $V=Y$.

By \cite[Lemma 3.9.3 (1)$\Rightarrow$(2)]{CC}
and by the assumption
that $Z$ is finite over $Y$,
the composition 
$g\circ f\colon X\to T$
has isolated characteristic points
at the inverse image $Z\times_Yv$.
Hence, 
the composition 
$g\circ f\colon X\to T$
is locally acyclic relatively to
${\cal F}$ on a neighborhood
of the fiber $X\times_Yv$
except at $Z\times_Yv$
and we have a canonical isomorphism
$$
\phi_v(Rf_*{\cal F},g)
\to
\bigoplus_{u\in Z\times_Yv}
\phi_u({\cal F},g\circ f).$$
Thus by the Milnor formula (\ref{eqMil}),
we have
\begin{align*}
-\dim{\rm tot} 
\phi_v(Rf_*{\cal F},g)
&=
\sum_{u\in Z\times_Yv}
-\dim{\rm tot} 
\phi_u({\cal F},g\circ f)
\\
&
=
\sum_{u\in Z\times_Yv}
(CC{\cal F},d(g\circ f))_{T^*X,u}
=
(f_*CC{\cal F},dg)_{T^*Y,v}
\end{align*}
and the assertion follows.
\qed}

\medskip

If $Y$ is a curve,
Conjecture \ref{cnf*}.2
may be rephrased as follows.
Let $C=SS{\cal F}$
be the singular support and
assume that
on a dense open subscheme
$V\subset Y$, the restriction $f_V\colon 
X_V=X\times_YV\to V$
of $f$ is $C$-transversal.
Then, 
$f_\circ C\times_YV$
is a subset of the $0$-section.
Thus the condition that
every irreducible component
of $f_\circ C$ is of dimension 1
is satisfied.
Further $f_V$ is locally acyclic
relatively to ${\cal F}$.
Since $f$ is proper, $Rf_*{\cal F}$
is locally constant on $V$
by \cite[Th\'eor\`eme 2.1]{Artin}
and we have
\begin{equation}
CCRf_*{\cal F}
=
-\Bigl(
{\rm rank}\ Rf_*{\cal F}\cdot
[T^*_YY]
+
\sum_{y\in Y\!\! \sm \!\! V}
a_yRf_*{\cal F}\cdot
[T^*_yY]
\Bigr).
\label{eqCCRf}
\end{equation}
For a closed point $y\in Y$,
the Artin conductor
$a_yRf_*{\cal F}$
is defined by
\begin{equation}
a_yRf_*{\cal F}
=
\chi(X_{\bar \eta},{\cal F})
-
\chi(X_{\bar y},{\cal F})
+
{\rm Sw}_y
H^*(X_{\bar \eta},{\cal F}).
\label{eqay}
\end{equation}
In the right hand side, the first two terms denote
the Euler-Poincar\'e characteristics
of the geometric generic fiber
and the geometric closed fiber
respectively
and the last term denotes
the Swan conductor at $y$.

Let $df$ denote the section of $T^*X$
on a neighborhood of 
the inverse image $X_y$
defined by the pull-back of
a basis $dt$ of the line bundle $T^*Y$ 
for a local coordinate $t$ on
a neighborhood of $y\in Y$.
Then, the intersection product
$(CC{\cal F},df)_{T^*X,y}$
supported on the inverse image of $X_y$
is well-defined since
$SS{\cal F}$ is a closed
conical subset.

\begin{lm}\label{lmproj}
Let $C=SS{\cal F}$
be the singular support and
let 
$V\subset Y$ be a dense open
subset such that
$f_V\colon X_V\to V$
is projective, smooth
and $C$-transversal.

{\rm 1.}
The equality
{\rm (\ref{eqcnf})} is equivalent
to the equality
\begin{equation}
-a_yRf_*{\cal F}
=
(CC{\cal F},df)_{T^*X,y}
\label{eqaf}
\end{equation}
at each point $y\in Y\sm V$,
where the right hand side denotes
the intersection number supported
on the inverse image of $y$.

{\rm 2.}
Further, if $f$ has at most 
isolated characteristic points,
then
Conjecture {\rm \ref{cnf*}.2} holds.
In particular,
if $f\colon X\to Y$
is a finite flat generically \'etale
morphism
of smooth curves, then 
Conjecture {\rm \ref{cnf*}.2} holds.

{\rm 3.}
Let $\delta_y$ denote
the difference of {\rm (\ref{eqaf})}.
If $X$ and $Y$ are projective,
we have $\sum_{y\in Y\!\!\sm \!\! V}\delta_y
\cdot \deg y=0$.
\end{lm}

\proof{
1.
Let $V'\subset V$ 
be the complement of
the images of irreducible
components $C_a$
of the singular support
$C=\bigcup_aC_a$
such that the image of
$C_a$ is a closed point of $Y$.
Then, for every closed point $w\in V'$,
the immersion $i_w\colon X_w\to X$
is properly $C$-transversal
and we have
$CCi_w^*{\cal F}=
i_w^!CC{\cal F}$
by \cite[Theorem 7.6]{CC}.
Further, we have
\begin{equation}
f_*CC{\cal F}
=
-\Bigl(
(i_w^!CC{\cal F},T^*_{X_w}X_w)
\cdot [T^*_YY]
-
\bigcup_{y\in Y\!\! \sm \!\! V}
(CC{\cal F},df)_{T^*X,y}
\cdot
[T^*_yY]
\Bigr).
\label{eqfCCV}
\end{equation}
If we assume that 
$f_V\colon X\times_YV\to V$
is projective,
the index formula (\ref{eqif})
implies 
$$
{\rm rank}\ Rf_*{\cal F}
=
\chi(X_w,i_w^*{\cal F})
=
(CCi_w^*{\cal F},T^*_{X_w}X_w)_{T^*X_w}.$$
Thus, it suffices
to compare
(\ref{eqfCCV})
and 
(\ref{eqCCRf}).

2.
If $f$ has at most 
isolated characteristic points,
(\ref{eqaf}) is an immediate consequence of
the Milnor formula (\ref{eqMil}).

3.
We have
\begin{align*}
(CCRf_*{\cal F},T^*_YY)_{T^*Y}
&\ =
\chi(Y,Rf_*{\cal F})\\
&\ =
\chi(X,{\cal F})
=
(CC{\cal F},T^*_XX)_{T^*X}
=
(f_*CC{\cal F},T^*_YY)_{T^*Y}
\end{align*}
by the index formula (\ref{eqif})
and the projection formula.
Thus it follows from 1.
\qed}
\medskip

We prove some cases of
Conjecture \ref{cnf*}.2 assuming
that $X$ is a surface.

Let $X$ be a normal noetherian
scheme and $U\subset X$ be a dense
open subscheme.
Let $G$ be a finite group and
$V\to U$ be a $G$-torsor.
The normalization $Y\to X$
in $V$ carries a natural action of
$G$.
For a geometric point $\bar x$
of $X$, 
the stabilizer $I\subset G$
of a geometric point $\bar y$
of $Y$ above $\bar x$
is called an inertia subgroup at $\bar x$.

\begin{lm}\label{lmSSj}
Let $G$ be a finite group and
\begin{equation}
\begin{CD}
W@<{j'}<<V\\
@VrVV@VV{r'}V\\
X@<j<<U
\end{CD}
\end{equation}
be a cartesian diagram
of smooth schemes over a field $k$
where the horizontal arrows
are dense open immersions,
the right vertical arrow
$V\to U$ is a $G$-torsor
and 
the left vertical arrow
$r\colon W\to X$
is proper.
Assume that for every geometric point
$x$ of $X$,
the order of the inertia group
$I_x\subset G$
is prime to $\ell$.

Let ${\cal F}$ be a locally constant
sheaf on $U$ such that
the pull-back $r^{\prime *}{\cal F}$
is a constant sheaf.
Then, for the intermediate
extension $j_{!*}{\cal F}
=j_{!*}({\cal F}[\dim U])[-\dim U]$ on $X$,
we have an inclusion
$$SSj_{!*}{\cal F}
\subset r_\circ (T^*_WW).$$
\end{lm}

\proof{
Since the assertion is \'etale local on
$X$, we may assume that
$G=I_x$ is of order prime to $\ell$.
Then, the canonical morphism
${\cal F}=
(r'_*r^{\prime*}{\cal F})^G
\to
r'_*r^{\prime*}{\cal F}$
is a splitting injection
and induces a splitting injection
$j_{!*}{\cal F}\to
j_{!*} r'_*r^{\prime*}{\cal F}$.
Hence, we have
$$SS(j_{!*}{\cal F})
\subset
SS(j_{!*} r'_*r^{\prime*}{\cal F}).$$
Since every irreducible
subquotient of the shifted perverse sheaf
$j_{!*} r'_*r^{\prime*}{\cal F}$ 
is isomorphic to an irreducible
subquotient of
a shifted perverse sheaf
$^p\!{\cal H}^0(Rr_*j'_*r^{\prime*}{\cal F}
[\dim U])[-\dim U]$
extending $r'_*r^{\prime*}{\cal F}$,
we have
$$SS(j_{!*} r'_*r^{\prime*}{\cal F})
\subset
SS(\ \! ^p\!{\cal H}^0(Rr_*j'_*r^{\prime*}{\cal F}))
\subset
SS(Rr'_*j'_*r^{\prime*}{\cal F})
$$
by \cite[Theorem 1.4 (ii)]{Be}.
Since 
$j'_*r^{\prime*}{\cal F}$
is a constant sheaf on $W$
and $r$ is proper,
we have
$$
SS(Rr_*j'_*r^{\prime*}{\cal F})
\subset
r_\circ SS(j'_*r^{\prime*}{\cal F})
\subset
r_\circ (T^*_WW)$$
by \cite[Lemma 2.2 (ii), Lemma 2.1 (iii)]{Be}.
Thus the assertion follows.
\qed}

\begin{pr}\label{pr1}
Let $X$ be a normal scheme
of finite type over a perfect field $k$,
$Y$ be a smooth curve
over $k$ and $f\colon X\to Y$
be a flat morphism over $k$.
Let $V\subset Y$ be a
dense open subscheme such
that $f_V\colon X_V=X\times_YV\to V$
is smooth.

{\rm 1.}
There exist a
finite flat surjective morphism
$g\colon Y'\to Y$ of smooth curves
over $k$
and a dense open subscheme
$X''\subset X'$ of the 
normalization $X'$ of
$X\times_YY'$
satisfying the following condition:

{\rm (1)}
We have inclusions
$X'_V=X'\times_YV
\subset X''\subset X'$
and $X''\subset X'$
is dense in every fiber of $X'\to Y'$.
The morphism $X''\to Y'$ is smooth.

{\rm 2.}
Let ${\cal F}$ be a perverse sheaf
on $X_V$ and let $C\subset T^*X_V$
be a closed conical subset
on which ${\cal F}$ is micro-supported.
Assume that 
$f_V\colon X_V\to V$
is $C$-transversal.
Then, there exist
$g\colon Y'\to Y$ 
and 
$X''\subset X'\to X\times_YY'$
as in {\rm 1.}\
satisfying the condition {\rm (1)}
above and the following condition:

{\rm (2)}
Let $j'\colon X'_V\to X''$
denote the open immersion
and $C'=SSj'_{!*}{\cal F}'$ 
be the singular support of 
the intermediate extension
$j'_{!*}{\cal F}'$ of the pull-back
${\cal F}'$ of ${\cal F}$ to
$X'_V$.
Then, 
the morphism
$X''\to Y'$ is $C'$-transversal.
\end{pr}

\proof
{By devissage and approximation,
we may assume that
the complement
$Y\sm V$ consists
of a single closed point $y$
and that the closed fiber
$X_y$ is irreducible.
The assertion is local
on a neighborhood in $X$ of
the generic point $\xi$
of $X_y$.

1. It follows from \cite{Epp}.

2.
Since $f_V\colon X_V\to V$ is smooth,
the $C$-transversality of $f_V$
and the condition that ${\cal F}$
is micro-supported on $C$
are preserved after base change
by \cite[Lemma 3.9.2, Lemma 4.2.4]{CC}.
After replacing $Y$ by $Y'$
and $X$ by $X''$ as in 1.,
we may assume 
that $X\to Y$ is smooth.
Shrinking $X$
and $Y$ further if necessary,
we may assume that
${\cal F}$ is locally constant.

Let $W_V\to X_V$ be a
$G$-torsor for a finite group $G$
such that the pull-back of
${\cal F}$ on $W_V$ is a constant sheaf.
Let $r\colon W\to X$ 
be the normalization
of $X$ in $W_V$.
Applying 1 to $W\to Y$
and shrinking $X$ if necessary,
we may assume that
there exists a finite flat surjective
morphism of smooth curves $Y'\to Y$
such that the normalization 
$W'$ of $W\times_YY'$ is smooth over $Y'$.

Let $r'\colon W'\to X'$
be the canonical morphism.
Since the ramification
index at the 
generic point $\xi'$ 
of an irreducible component
of the fiber $X'_y$
is $1$,
the inertia group at $\xi'$
is of order a power of $p$.
Hence, after shrinking $X$ if necessary,
we may assume that for every geometric point
$w'$ of $W'$, the order of
the inertia group is a power of $p$.
Hence, by Lemma \ref{lmSSj},
we have 
$C'=SS(j'_{!*}{\cal F})
\subset r'_\circ  (T^*_{W'}W')$.
Since $W'\to Y'$ is smooth,
the morphism $f'\colon X'\to Y'$
is $C'$-transversal by
\cite[Lemma 3.9.3 (2)$\Rightarrow$(1)]{CC}.
\qed}

\begin{thm}\label{thm1}
Let the notation be as in Conjecture
{\rm \ref{cnf*}}
and let $C=SS{\cal F}$ be
the singular support.
Assume that $\dim X=2,\ \dim Y=1$
and that there exists
a dense open subscheme $V\subset Y$
such that $f_V\colon X_V\to V$
is smooth and $C$-transversal.
Then, Conjecture {\rm \ref{cnf*}.2}
holds.
\end{thm}

\proof{
We may assume ${\cal F}$ is a
perverse sheaf by \cite[Theorem 1.4 (ii)]{Be}.
Since the resolution of singularity
is known for curves and surfaces,
we may assume $Y$ is projective.
Since a proper smooth surface
over a field is projective,
the surface $X$ is projective.
Let $y\in Y\sm V$
be a point. It suffices to show
the equality (\ref{eqaf}).

By Proposition \ref{pr1}
and approximation,
there exists a finite flat surjective
morphism $Y'\to Y$ of
proper smooth curves
\'etale at $y$
and satisfying the conditions in 
Proposition \ref{pr1} on the
complement $Y\sm \{y\}$.
Since the normalization 
$X'$ of $X\times_YY'$
is projective,
we may take a projective smooth scheme
$P$ and decompose $f'\colon X'\to Y$
as a composition $X'\to P\to Y$
of a closed immersion $i\colon X'\to P$
and $g\colon P\to Y$.

Let $U=V\cup\{y\}$.
Let ${\cal F}'$ be
the pull-back of ${\cal F}$
to $X'_U=X'\times_YU$
and let $j'_{!*}{\cal F}'$
be the intermediate extension
with respect to the open immersion
$j'\colon X'_U\to X'$.
It suffices to show
Conjecture \ref{cnf*}.2 holds
for $P\to Y$ and ${\cal G}=i_*j'_{!*}{\cal F}'$.
Outside the inverse image
of $y$,
the morphism $g\colon P\to Y$
has at most isolated characteristic points
with respect to the singular support
$SS{\cal G}$
by the condition (2) in Proposition \ref{pr1}
and 
\cite[Lemma 3.9.3 (2)$\Rightarrow$(1)]{CC}
applied to the restriction of
the immersion $X'\to P$
on the complement of a finite closed
subset of $X'$.
Thus, we have $\delta_{y'}=0$
for any closed point $y'\in Y'$ not on $y$
by Lemma \ref{lmproj}.2.
This implies
$[Y'\colon Y]\cdot \delta_y=0$
by Lemma \ref{lmproj}.3.
Thus 
the assertion follows.
\qed}

\begin{thm}\label{thm2}
Let the notation be as in Conjecture
{\rm \ref{cnf*}}
and let $C=SS{\cal F}$ be
the singular support.
Assume that $\dim X=\dim Y=2$,
that $f\colon X\to Y$ is proper surjective
and that every irreducible component
of $f_\circ C$ is of dimension $2$.
Then, Conjecture {\rm \ref{cnf*}.2}
holds.
\end{thm}

\proof{
By Lemma \ref{lmfini},
the assertion holds
except possibly for
the coefficients of
the fibers $T^*_yY$
of finitely many closed points
$y\in Y$ where
$X\to Y$ is not finite.
Let $v\in Y$ be a closed point
and we show that
the coefficients of
the fibers $T^*_vY$ are equal.

Since the resolution of singularity
is known for surfaces
and since a proper smooth surface
over a field is projective,
we may assume that $Y$
and hence $X$ are projective.
By replacing $X$ by the Stein factorization
of $X\to Y$ except on
a neighborhood of $v$,
we define $X\to X'\to Y$
such that $f'\colon X'\to Y$ is finite
on the complement of $v$
and $r\colon X\to X'$ is an isomorphism
on the inverse image of
a neighborhood of $v$.

Since $X'$ is projective,
we may take a projective smooth scheme
$P$ and decompose $f'\colon X'\to Y$
as a composition $X'\to P\to Y$
of a closed immersion $i\colon X'\to P$
and $g\colon P\to Y$.
Conjecture \ref{cnf*}.2 holds
for ${\cal G}=i_*Rr_*{\cal F}$ and
$g\colon P\to Y$
except possibly for
the coefficients of
the fiber $T^*_vY$
by Lemma \ref{lmfini}.
Namely, we have
$g_*CC{\cal G}
=CCRg_*{\cal G}
=CCRf_*{\cal F}$
except possibly for
the coefficients of
the fiber $T^*_vY$.
By the index formula (\ref{eqif}),
we have
\begin{align*}
(CCRf_*{\cal F},T^*_YY)_{T^*Y}
&\ =
\chi(Y,Rf_*{\cal F})\\
&\ =
\chi(P,{\cal G})=
(CC{\cal G},T^*_PP)_{T^*P}
=
(g_*CC{\cal G},T^*_YY)_{T^*Y}.
\end{align*}
Thus, 
we have an equality
also for the coefficients of
the fiber $T^*_vY$.
\qed}

\medskip

We give a characterization of
characteristic cycle
using functoriality.
For the definition of
$C$-transversal morphisms
$h\colon W\to X$
of smooth schemes
and 
the pull-back
$Z_n(C)\to Z_m(h^\circ C)$,
we refer to \cite[Definition 7.1]{CC}.
For a constructible complex ${\cal F}$
on a projective space ${\mathbf P}={\mathbf P}^n$,
let $R {\cal F}
=Rp^\vee_*p^*{\cal F}[n-1]$
denote the Radon transform 
on the dual projective space
${\mathbf P}^\vee$
where $p\colon Q\to {\mathbf P}$
and $p^\vee\colon Q\to {\mathbf P}^\vee$
denote the projections
on the universal family of
hyperplanes $Q=\{(x,H)\in 
{\mathbf P}\times {\mathbf P}^\vee
\mid x\in H\}$.
For a linear combination $A=\sum_am_aC_a$
of irreducible closed conical
subset $C_a\subset T^*{\mathbf P}$
of dimension $n$,
let $L A=(-1)^{n-1}p^\vee_*p^!A$
denote the Legendre transform
(cf.\! \cite[Corollary 7.5]{CC}).

\begin{pr}\label{prch}
Let $k$ be a perfect field
and $\Lambda$ be a finite
field of characteristic $\ell$ invertible in $k$.
Then, there exists a unique way
to attach a linear combination
$A({\cal F})=\sum_am_aC_a$
satisfying the conditions {\rm (1)-(5)} below
of irreducible components
of the singular support
$SS{\cal F}=C=\bigcup_aC_a\subset T^*X$
to each
smooth scheme $X$ over $k$
and each constructible complex ${\cal F}$
of $\Lambda$-modules on $X$:

{\rm (1)}
For every \'etale morphism
$j\colon U\to X$,
we have $A(j^*{\cal F})=
j^*A({\cal F})$.

{\rm (2)}
For every properly $C$-transversal
closed immersion
$i\colon W\to X$
of smooth schemes,
we have
$A(i^*{\cal F})=
i^!A({\cal F})$.

{\rm (3)}
For every closed immersion
$i\colon X\to P$
of smooth schemes,
we have
$A(i_*{\cal F})=
i_*A({\cal F})$.

{\rm (4)}
For the Radon transform,
we have $A(R {\cal F})=
LA({\cal F})$.

{\rm (5)}
For $X={\rm Spec}\ k$,
we have
$A({\cal F})={\rm rank}\ {\cal F}\cdot T^*_XX$.

\noindent
If $A({\cal F})$
satisfies the conditions {\rm (1)--(5)},
then we have
\begin{equation}
A({\cal F})=CC{\cal F}.
\label{eqA=C}
\end{equation}
\end{pr}

As the proof below shows,
it suffices to assume the condition (2)
in the case where 
$i\colon x\to X$ is the closed immersion
of a closed point or
$i\colon L\to {\mathbf P}$ 
is the closed immersions of lines
in projective spaces.
By the definition of
the naive Radon transform
and the naive Legendre transform,
the equality in 
condition (4) can be decomposed as
$A(Rp^\vee_*p^*{\cal F})=
p^\vee_*A(p^*{\cal F})=
p^\vee_*p^! A({\cal F})$.
The first (resp.\! second) equality 
corresponds to a special case 
of Conjecture \ref{cnf*}.2
by \cite[Lemma 3.11]{CC}
(resp.\! to a special case of
\cite[Proposition 5.17]{CC}).

\proof{
The characteristic cycles
satisfy the conditions (1)--(5)
by
\cite[Lemma 5.11.2]{CC},
\cite[Theorem 6.6]{CC},
\cite[Lemma 5.13.2]{CC},
\cite[Corollary 7.12]{CC},
\cite[Lemma 5.11.1]{CC}
respectively.
Thus the existence is proved.

We show the uniqueness.
It suffices to show the equality (\ref{eqA=C}).
By the condition (2)
applied to the closed immersion
$x\to X$
of a closed point in a dense open 
subscheme $U\subset X$
where ${\cal F}$ is locally constant
and by the conditions (1) and (5),
the coefficient of the $0$-section $T^*_XX$
in $A({\cal F})$
equals the rank of the restriction
${\cal F}|_U$.

We show that this property
and the conditions
(3)--(5) imply the index formula
for projective smooth scheme $X$.
By (3), we may assume
that $X$ is a projective space ${\mathbf P}
={\mathbf P}^n$ for $n\geqq 2$.
Let $R^\vee {\cal G}=
Rp_*p^{\vee*}{\cal G}(n-1)[n-1]$
and $L^\vee B=(-1)^{n-1}p_*p^{\vee!} B$
denote the inverse Radon transform 
and  the inverse Legendre transform.
Then, 
$R^\vee R {\cal F}$
is isomorphic to ${\cal F}$ 
up to locally constant complex of
rank 
$(n-2)\cdot \chi({\mathbf P}_{\bar k},{\cal F})$.
Hence the coefficient of
the $0$-section $T^*_{\mathbf P}{\mathbf P}$
in $A(R^\vee R {\cal F})$
equals that in $A({\cal F})$
plus $(n-1)\cdot \chi({\mathbf P}_{\bar k},{\cal F})$.
Similarly, the coefficient of
$T^*_{\mathbf P}{\mathbf P}$
in $L^\vee LA ({\cal F})$
equals that in $A({\cal F})$
plus $(n-1)\cdot (A({\cal F}), T^*_{\mathbf P}{\mathbf P})
_{T^*{\mathbf P}}$.
Thus, we obtain the index formula
$\chi({\mathbf P}_{\bar k},{\cal F})
=(A({\cal F}), T^*_{\mathbf P}{\mathbf P})
_{T^*{\mathbf P}}$
for projective $X$.

Next, we show
(\ref{eqA=C}) assuming $\dim X=1$.
Let $U\subset X$ be
a dense open subscheme 
where ${\cal F}$ 
is locally constant.
Then, we have
$$A({\cal F})
=-\Bigl({\rm rank}\ {\cal F}
\cdot T^*_XX+
\sum_{x\in X\!\!\sm \!\! U}
m_x\cdot T^*_xX\Bigr).
$$
We show $m_x$
equals the Artin
conductor $a_x{\cal F}$.
If ${\cal F}|_U$ is unramified at $x$
and if ${\cal F}_x=0$,
the same argument as above shows
$m_x={\rm rank}\ {\cal F}|_U=a_x{\cal F}$.
We show the general case.
By (1), we may assume $U=X\sm \{x\}$.
Further shrinking $X$,
we may assume that there
exists a finite \'etale surjective morphism
$\pi\colon X'\to X$ such that
the pull-back
$\pi^*{\cal F}$ is unramified at 
every point $\bar X'\sm X'$
of the boundary of a smooth
compactification $j'\colon X'\to \bar X'$.
Then, we have
$$A(j'_*\pi^*{\cal F})
=-
\Bigl({\rm rank}\ {\cal F}|_U\cdot
T^*_{\bar X'}\bar X'+
\sum_{x'\in \pi^{-1}(x)}
m_x\cdot 
T^*_{x'}\bar X'
+
\sum_{x'\in \bar X'\sm X'}
a_x{\cal F}\cdot
T^*_{x'}\bar X'
\Bigr).$$
Thus, the index formula implies
$m_x=a_x{\cal F}$.

We show the general case.
By (1), we may assume $X$
is affine.
We consider an immersion
$X\to {\mathbf A}^n\subset {\mathbf P}^n$.
Then, by (1) and (3),
we may assume $X$ is projective.
Set $SS{\cal F}=
C=\bigcup_aC_a$ and
take a projective embedding $X\to 
{\mathbf P}$ and
a pencil $L\subset {\mathbf P}^\vee$
satisfying the following properties
as in \cite[Lemma 2.3]{SY}:
The axis $A_L$ of the pencil meets $X$
transversely, that
the blow-up $\pi_L\colon X_L\to X$
is $C$-transversal,
that the morphism
$p_L\colon X_L\to L$
defined by the pencil has at most
isolated characteristic points,
that the isolated characteristic points
are not contained in the inverse image
of $V$ and are unique in the fibers
of $p_L$,
and that for each irreducible component
$C_a$ there exists an
isolated characteristic point $u$
where a section $dp_L$ of $T^*X$ meets
$C_a$.

Then, for $v=p_L(u)$,
the coefficient $m_v$ of $T^*_vL$
in $A(Rp_{L*}\pi^*{\cal F})$
equals the Artin conductor
$-\dim{\rm tot}\phi_u({\cal F},p_L)=
-a_vRp_{L*}\pi^*{\cal F}.$
Let $i_L\colon L\to {\mathbf P}^\vee$
denote the immersion.
Then,
by the proper base change theorem
and by the conditions (2) and (3),
we have
$A(Rp_{L*}\pi^*{\cal F})
=(-1)^{n-1}i_L^!A(R{\cal F})
=i_L^!L A({\cal F})$.
Thus, if
$A({\cal F})=\sum_am_aC_a$,
the coefficient $m_v$ equals
$m_a\cdot (C_a,dp_L)_{T^*X,u}
\neq 0$
and we obtain
$$-\dim{\rm tot}\phi_u({\cal F},p_L)=
m_a\cdot (C_a,dp_L)_{T^*X,u}
\neq 0.$$
This means that the coefficient
$m_a$ is characterized by
the same condition as the
Milnor formula (\ref{eqMil})
and we have (\ref{eqA=C}).
\qed}

\end{document}